\documentclass[twocolumn]{autart}    % Enable this line and disable the
\usepackage{amssymb,latexsym,color,amsmath,pifont,epsfig,graphicx}
\usepackage{setspace,cite}
\usepackage{algorithm}
\usepackage{algorithmic}
\usepackage{subfig}  %Use as needed
\usepackage{algorithm}
\usepackage{algorithmic}
\usepackage{soul}
\usepackage{graphics}
\usepackage{multirow}
\usepackage{bm}
\usepackage{rotating}
\usepackage{epsfig}
\usepackage{amssymb,latexsym,color,amsmath,pifont,epsfig,mathtools}
\usepackage{indentfirst}
\usepackage{dsfont}
\usepackage[justification=centering]{caption}
\usepackage{url}
\usepackage{tikz}
\usepackage{graphicx,booktabs,ctable}
\usetikzlibrary{arrows,shapes,trees,calc,positioning}
\usetikzlibrary{matrix}
\DeclareMathOperator*{\argmin}{argmin}

\DeclareMathOperator*{\minimize}{minimize}

\DeclareMathOperator*{\subject}{subject~to}

\newcommand{\DefinedAs}[0]{\mathrel{\mathop:}=}
\newcommand{\bbR}{\mathbb{R}}
\newcommand{\tc}{\textcolor}

\newcommand{\mre}{\mathrm{e}}
\newcommand{\mrd}{\mathrm{d}}

\newcommand{\prox}{\mathbf{prox}}

\newtheorem{theorem}{Theorem}

\newtheorem{lemma}[theorem]{Lemma}

\newtheorem{remark}{Remark}

\newcommand{\enma}[1]   {\ensuremath{#1}}

\newcommand{\non}{\nonumber}

\newcommand{\beq}{\begin{equation}}
\newcommand{\eeq}{\end{equation}}
\newcommand{\bseq}{\begin{subequations}}
\newcommand{\eseq}{\end{subequations}}
\newcommand{\beqn}{\begin{eqnarray}}
\newcommand{\eeqn}{\end{eqnarray}}
\newcommand{\ba}{\begin{array}}
\newcommand{\ea}{\end{array}}
\newcommand{\bct}{\begin{center}}
\newcommand{\ect}{\end{center}}
\newcommand{\btmz}{\begin{itemize}}
\newcommand{\etmz}{\end{itemize}}
\newcommand{\benum}{\begin{enumerate}}
\newcommand{\eenum}{\end{enumerate}}

\newcommand{\cL}{\enma{\mathcal L}}

\newcommand{\norm}[1]{\| #1 \|}                 %does not make large \|

\newcommand{\inner}[2]{\left\langle #1,#2 \right\rangle}

\newcommand{\matbegin}{
        \left[
}
\newcommand{\matend}{
        \right]
}

\newcommand{\tbo}[2]{
  \matbegin \begin{array}{c}
       #1 \\ #2
       \end{array} \matend }

\newcommand{\tbt}[4]{
  \matbegin \begin{array}{cc}
       #1 & #2 \\ #3 & #4
       \end{array} \matend }

\newcommand{\be}{\begin{equation}}
\newcommand{\ee}{\end{equation}}

\newcommand{\cplxs}{ C\kern -.35em \rule{0.03 em}{.7 ex}~   }

\def\complex{\hbox{C\kern -.45em \rule{0.03 em}{1.5 ex}}~}

\newcommand{\bi}{\begin{itemize}}
\newcommand{\ei}{\end{itemize}}

\newtheorem{assumption}{Assumption}

\begin{document}
%\tableofcontents

\begin{frontmatter}
%\runtitle{Insert a suggested running title}  % Running title for regular
                                              % papers but only if the title
                                              % is over 5 words. Running title
                                              % is not shown in output.

\title{\large Proximal gradient flow and Douglas-Rachford splitting dynamics:
	\\[0.1cm]
global exponential stability via integral quadratic constraints\thanksref{footnoteinfo}}

\thanks[footnoteinfo]{Supported in part by the National Science Foundation under awards ECCS-1708906 and ECCS-1809833.}

	\vspace*{-0.75cm}
\author{Sepideh Hassan-Moghaddam}\ead{hassanmo@usc.edu}~and
\author[mj]{Mihailo R.\ Jovanovi\'c}\ead{mihailo@usc.edu}

\address[mj]{Ming Hsieh Department of Electrical Engineering, University of Southern California, Los Angeles, CA 90089}

\begin{abstract}
Many large-scale and distributed optimization problems can be brought into a composite form in which the objective function is given by the sum of a smooth term and a nonsmooth regularizer. Such problems can be solved via a proximal gradient method and its variants, thereby generalizing gradient descent to a nonsmooth setup. In this paper, we view proximal algorithms as dynamical systems and leverage techniques from control theory to study their global properties. In particular, for problems with strongly convex objective functions, we utilize the theory of integral quadratic constraints to prove the global exponential stability of \tc{black}{the equilibrium points of} the differential equations that govern the evolution of proximal gradient and Douglas-Rachford splitting flows. In our analysis, we use the fact that these algorithms can be interpreted as variable-metric gradient methods on \tc{black}{the suitable} envelopes and exploit structural properties of the nonlinear terms that arise from the gradient of the smooth part of the objective function and the proximal operator associated with the nonsmooth regularizer. We also demonstrate that these envelopes can be obtained from the augmented Lagrangian associated with the original nonsmooth problem and establish conditions for global exponential convergence even in the absence of strong convexity.
\end{abstract}

    \vspace*{-0.5cm}

    \end{frontmatter}

	\vspace*{-8ex}
\section{Introduction}
\label{sec.intro}

	\vspace*{-2ex}
Structured optimal control and estimation problems typically lead to optimization of objective functions that consist of a sum of a smooth term and a nonsmooth regularizer. Such problems are of increasing importance in applications and it is thus necessary to develop efficient algorithms for distributed and embedded nonsmooth composite optimization~\cite{nedozd09,waneli11,latstepat16,latfrepat19}. The lack of differentiability in the objective function precludes the use of standard descent methods from smooth optimization. Proximal gradient method~\cite{becteb09,parboy13} generalizes gradient descent to nonsmooth context and provides a powerful tool for solving problems in which the nonsmooth term is separable over the optimization variable.

Examining optimization algorithms as continuous-time dynamical systems has been an active topic since the seminal work of Arrow, Hurwicz, and Uzawa~\cite{arrhuruza58}. This viewpoint can provide important insight into performance of optimization algorithms and \tc{black}{streamline their convergence analysis. During}~the last decade, it has been advanced and extended to a broad class of problems including convergence analysis of primal-dual~\cite{feipagAUT10,waneli11,chemalcorSCL16,chemallowcor18,dhikhojovTAC19,quli18} and accelerated~\cite{suboycan16,wibwiljor16,frarobvid18,shidujorsu18,muejor19,povli19} first-order methods. Furthermore, establishing the connection between theory of ordinary differential equations (ODEs) and numerical optimization algorithms has been a topic of many studies, including~\cite{brobar89,schsin00}; for recent efforts, see~\cite{wibwiljor16,zhamoksrajad18}. 

Optimization algorithms can be viewed as a feedback interconnection of linear dynamical systems with nonlinearities that possess certain structural properties. This system-theoretic interpretation was exploited in~\cite{lesrecpac16} and further advanced in recent papers~\cite{dhikhojovTAC19,huseiran17,hules17,fazribmorpre18,mogjovACC18,mogjovCDC18a,dinhudhijovCDC18,seifazprepap19}. The key idea is to exploit structural features of linear and nonlinear terms and utilize theory and techniques from stability analysis of nonlinear dynamical systems to study properties of optimization algorithms. This approach provides new methods for studying not only convergence rate but also robustness of optimization routines~\cite{mohrazjovCDC18,mohrazjovACC19,mohrazjovTAC19,micschebe19} and can lead to new classes of algorithms that strike a desired tradeoff between the speed and robustness.

In this paper, we utilize techniques from control theory to establish global properties of proximal gradient flow and Douglas-Rachford (DR) splitting dynamics. These algorithms provide an effective tool for solving nonsmooth convex optimization problems in which the objective function is given by a sum of a differentiable term and a nondifferentiable regularizer. When \tc{black}{the} smooth \tc{black}{term} is strongly convex with a Lipschitz continuous gradient, we prove the global exponential stability \tc{black}{of both the proximal gradient flow and the DR splitting dynamics} by utilizing the theory of IQCs~\cite{megran97}. We also generalize the Polyak-Lojasiewicz (PL)~\cite{pol63} condition to nonsmooth problems and show global exponential convergence of the \tc{black}{forward-backward (FB) envelope~\cite{patstebem14,stethepat17,thestepat18}} even in the absence of strong convexity. 

\tc{black}{Although there are related approaches for studying optimization algorithms from a control-theoretic perspective, to the best of our knowledge, we are the first to introduce the continuous forms of proximal gradient and DR splitting algorithms. We use simple proofs to establish their global stability properties and provide explicit bounds on convergence rates. Furthermore, standard forms of these algorithms are obtained via explicit forward Euler discretization of \mbox{continuous-time dynamics.}}

The paper is structured as follows. In Section~\ref{sec.setup}, we formulate the nonsmooth composite optimization problem and provide background material. In Section~\ref{sec.proxalgs}, we establish the global exponential stability of the proximal gradient flow dynamics for a problem with strongly convex objective function. Moreover, by exploiting the problem structure, we demonstrate the global exponential convergence of the forward-backward envelope even in the absence of strong convexity. In Section~\ref{section.DRS}, we introduce a continuous-time gradient flow dynamics based on the celebrated Douglas-Rachford splitting algorithm and utilize the theory of IQCs to prove global exponential stability for strongly convex problems. We offer concluding remarks in Section~\ref{sec.conclusion}.

	\vspace*{-2ex}
\section{Problem formulation and background}
\label{sec.setup}

	\vspace*{-2ex}
We consider a composite optimization problem,
\beq
    \ba{rcl}
	&\minimize\limits_{x}&
	f(x)
	\; + \;
	g(T x)
    \ea
    \label{eq.main}
	\eeq
where $x \in \bbR^n$ is the optimization variable, $T \in \bbR^{m \times n}$ is a given matrix, $f$: $\bbR^n \to \bbR$ is a \tc{black}{convex} function with a Lipschitz continuous gradient, and $g$: $\bbR^m \to \bbR$ is a nondifferentiable convex function. Such optimization problems arise in a number of applications and depending on the structure of the functions $f$ and $g$, different first- and second-order algorithms can be employed to solve them. We are interested in studying global convergence properties of methods based on proximal gradient flow algorithms. In what follows, we provide background material that we utilize in the rest of the paper.

	\vspace*{-2ex}
\subsection{Proximal operator and the associated envelopes}

	\vspace*{-2ex}
The proximal operator of a proper, closed, and convex function $g$ is defined as
	\be
    \prox_{\mu g}(v)
    \;\DefinedAs\;
    \argmin \limits_z\,
    \left(g(z)
    \; + \;
    \dfrac{1}{2\mu}
    \,
    \|z\,-\,v\|_2^2\right)
    \label{eq.prox}
    \ee
\tc{black}{where $\mu$ is a positive parameter and $v$ is a given vector. It is determined by the resolvent operator associated with $\mu \partial g$, $\prox_{\mu g} \DefinedAs (I + \mu \partial g)^{-1}$, and is a single-valued firmly non-expansive mapping~\cite{parboy13}, i.e., for any $u$ and $v$,}
    \be
    \tc{black}{
    \ba{l}
    \| \prox_{\mu g}(u)\,-\,\prox_{\mu g}(v)\|^{2}_2
    \,\leq \,
    \\[.cm]
    \;\;\;\;\;\;\;\;\;\;\;\;\;\;\;\;\;\;\;\;\;\;\;\;\;\;\;\;\;\;\;\;\;\inner{ u\,-\,v}{\prox_{\mu g}(u)\,-\,\prox_{\mu g}(v)}.
    \ea}
    \non
    \ee
The value function of the optimization problem~\eqref{eq.prox} determines the associated Moreau envelope,
    \be
    M_{\mu g}(v)
    \;\DefinedAs\;
    g(\prox_{\mu g}(v))
    \,+ \,
    \dfrac{1}{2\mu}\, \|\prox_{\mu g}(v)\,-\,v\|_2^2
    \non
    \ee
which is a continuously differentiable function even when $g$ is not~\cite{parboy13}, with
	 $
    \mu \nabla M_{\mu g}(v)
     =
    v - \prox_{\mu g}(v).
    $

By introducing an auxiliary optimization variable $z$, problem~\eqref{eq.main} can be rewritten as follows,
    \beq
	\ba{rl}
	\minimize\limits_{x,\,z}
	&
	f(x)
	\; + \;
	g(z)
    \\[.cm]
    \subject
    &
    T x\,-\,z
    \;=\;
    0
	\ea
    \label{eq.main1}
	\eeq
and the associated augmented Lagrangian is given by,
	 \be
         \mathcal{L}_{\mu}(x, z; y)
    \, \DefinedAs \,
    f(x)
    \, + \,
    g(z)
    \, + \,
    \inner{y}{T x -z}
    \, + \,
    \tfrac{1}{2\mu}
    \,
    \|T x - z\|_2^2.
    \non
    \ee
The completion of squares yields,
    \be
    \cL_\mu
    \, = \,
    f(x)
    \, + \,
    g(z)
    \, + \,
    \tfrac{1}{2\mu}
    \,
    \norm{z  \, - \, (T x \, + \, \mu y)}_2^2
     \, - \,
    \tfrac{\mu}{2} \, \norm{y}_2^2
   	\non
    \ee
where $y$ is the Lagrange multiplier. The minimizer of $\cL_\mu$ with respect to $z$ is
	\be
	z^\star(x,y)
	\; = \;
	\prox_{\mu g}(T x \, + \, \mu y)
	%\label{eq.zstar}
	\non
	\ee
and the evaluation of $\cL_\mu$ along the manifold resulting from this explicit minimization yields the proximal augmented Lagrangian~\cite{dhikhojovTAC19}, 
	$
	\cL_\mu (x; y)
	\DefinedAs 
	\cL_\mu(x, z^\star(x,y); y),
	$
	\be
	\cL_\mu (x; y)
	\; = \;
	f(x)
    	\; + \;
    	M_{\mu g}
    	(T x \, + \, \mu y)
    	\; - \;
    	\tfrac{\mu}{2} \, \| y \|^2_2.
	\label{eq.pal}
	\ee
This function is continuously differentiable with respect to both $x$ and $y$ and it can be used as a foundation for the development of first- and second-order primal-dual methods for nonsmooth composite optimization~\cite{dhikhojovTAC19,dhikhojovTAC17}. For $T = I$, the forward-backward envelope~\cite{patstebem14,stethepat17,thestepat18} is obtained by restricting the proximal augmented Lagrangian $\cL_\mu (x; y)$ along the manifold $y^\star (x) = - \nabla f(x)$ resulting from the KKT optimality conditions,
	\be
	\ba{rrl}
	F_\mu (x)
	&\! \DefinedAs \! &
	\cL_\mu(x; y^\star(x))
	\; = \;
	\cL_\mu(x; y = - \nabla f(x))
	\\[0.cm]
	&\! = \! &
	f(x)
    	\; + \;
    	M_{\mu g}
    	(x \, - \, \mu \nabla f(x))
    	\; - \;
    	\tfrac{\mu}{2} \, \| \nabla f(x) \|^2_2.
	\ea
	\non
 	\ee
		
	\vspace*{-3ex}		
\subsection{Strong convexity and Lipschitz continuity}

	\vspace*{-2ex}		
The function $f$ is $m_f$-strongly convex if
	 \be
    f(\hat{x})
    \; \geq \;
    f(x)
    \; + \;
    \inner{\nabla f(x)}{\hat{x}\,-\,x}
    \; + \;
    \dfrac{m_f}{2}\,\|\hat{x}\,-\,x\|_2^2
    \non
    \ee   
\tc{black}{and its gradient is $L_f$-Lipschitz continuous if} 
    \be
    \tc{black}{f(\hat{x})
    \; \leq \;
    f(x)
    \; + \;
    \inner{\nabla f(x)}{\hat{x}\,-\,x}
    \; + \;
    \dfrac{L_f}{2}\,\|\hat{x}\,-\,x\|_2^2}
    \non
    \ee
for any $x$ and $\hat{x}$. \tc{black}{When both properties hold we have}
    \be
    m_f \| x\,-\,\hat{x} \|_2
    \,\leq\,
    \| \nabla f(x) -\nabla f(\hat{x}) \|_2
    \,\leq\,
    L_f \|x\,-\,\hat{x}\|_2.
    \label{eq.Lf-mf}
    \ee
and the following inequality is satisfied~\cite{nes13},
    \be
    \ba{rcl}
    \!\inner{\nabla f(x) -\nabla f(\hat{x})}{x\,-\,\hat{x}}
    &\! \geq \! &
    \dfrac{m_f\,L_f}{m_f+L_f}\,\|x\,-\,\hat{x}\|_2^2
    ~ +
    \\[.25cm]
    & \!\! &
    \dfrac{1}{m_f+L_f}\,\|\nabla f(x)-\nabla f(\hat{x})\|_2^2.
    \ea
    \label{eq.Lfmf}
    \ee
Furthermore, the subgradient $\partial g$ of a nondifferentiable function $g$ is defined as the set of points $z \in \partial g(x)$ that for any $x$ and $\hat{x}$ satisfy,
    \be
    g(\hat{x})
    \, \geq \,
    g(x)
    \, + \,
    z^T (\hat{x}\,-\,x).
    \label{eq.subgrad}
    \ee

	\vspace*{-4ex}
\subsection{Proximal Polyak-Lojasiewicz inequality}
\label{sec.PL}

		\vspace*{-2ex}
The Polyak-Lojasiewicz (PL) condition can be used to prove linear convergence of a gradient descent \tc{black}{even in the absence of convexity}~\cite{karnutsch16}. For an unconstrained optimization problem with a non-empty solution set \tc{black}{and a twice differentiable objective function $f$ with a Lipschitz continuous gradient,} the PL condition is given by
    \be
    \|\nabla f(x)\|_2^2
    \; \geq \;
    \gamma\,( f(x)\,-\,f^\star )
    \non
    \ee
where $\gamma>0$ and $f^{\star}$ is the optimal value of~$f$. For nonsmooth optimization problem~\eqref{eq.main1} with $T = I$, the proximal PL inequality holds for $\mu \in (0,1/L_f)$ if there exist $\gamma > 0$ such that
    \be
    \|G_{\mu}(x)\|_2^2
    \,\geq\,
    \gamma\, (F_{\mu} (x)\,-\, F_{\mu}^\star).
    \label{eq.lbonG}
    \ee
Here, $L_f$ is the Lipschitz constant of $\nabla f$, $F_{\mu}$ is the FB envelope, and $G_\mu$ is the generalized gradient map,
	\be
    \tc{black}{G_{\mu}(x)
    \;\DefinedAs\;
    \dfrac{1}{\mu}
    \,
    (
    x\,-\,\prox_{\mu g}(x\,-\, \mu \nabla f(x))
    ).}
    \label{eq.G}
    \ee
\tc{black}{When $f$ is twice continuously differentiable, the FB envelope $F_\mu$ is continuously differentiable with~\cite{patstebem14},
    \be
    \nabla F_{\mu}(x)
    \;=\;
    (I\,-\,\mu \nabla^2 f(x))\,G_{\mu}(x).
    \label{eq.gradFB}
    \ee}
    
	\vspace*{-2ex}
\section{Exponential stability of proximal algorithms}
\label{sec.proxalgs}

		\vspace*{-2ex}
In this section, we briefly discuss the Arrow-Hurwicz-Uzawa gradient flow dynamics that can be used to solve~\eqref{eq.main1} by computing the saddle points of the proximal augmented Lagrangian~\cite{dhikhojovTAC19}. We then show that the proximal gradient method in continuous time can be  obtained from the proximal augmented Lagrangian method by restricting the dual variable along the manifold $y = - \nabla f(x)$. Finally, we discuss global stability properties of proximal algorithms both in the presence and in the absence of strong convexity.

Continuous differentiability of the proximal augmented Lagrangian~\eqref{eq.pal} can be utilized to compute its saddle points via the Arrow-Hurwicz-Uzawa dynamic,
	\be
	\left[
    \ba{c}
    \dot{x}
    \\[-0.1cm]
    \dot{y}
    \ea
    \right]
    \, = \,
    \left[
    \ba{c}
    -\mu \, (\nabla f(x) \;+\; T^T\nabla M_{\mu g}(Tx \,+\, \mu y) )
    \\[-0.1cm]
    \mu \, ( \nabla M_{\mu g}(Tx \,+\, \mu y) \;-\; y )
     \ea
    \right].
    \label{eq.GF}
    \ee
As shown in~\cite{dhikhojovTAC19}, \tc{black}{the optimal primal-dual pair ($x^\star,y^\star$) is the globally exponentially stable equilibrium point of~\eqref{eq.GF} and $x^\star$ is the solution of~\eqref{eq.main}} for convex problems in which the matrix $T T^T$ is invertible and the smooth part of the objective function $f$ is strongly convex.

For convex problems with $T = I$ in~\eqref{eq.main},
	\beq
	\minimize\limits_{x}
	~
	f(x)
	\; + \;
	g(x)
    \label{eq.main11}
	\eeq
\tc{black}{the optimality condition is given by
   \be
    0
    \; \in \;
    \nabla f(x^\star) \; + \; \partial g(x^\star)
    \label{eq.OC}
    \ee
Multiplying by $\mu$ and adding/subtracting $x^\star$ yields,
	 \[
    0
    \; \in \;
    \left[
    I \, + \, \mu \partial g
    \right] (x^\star)
    \; + \;
    \mu \nabla f(x^\star) \; - \; x^\star.
    \]
Since $\prox_{\mu g}$ is determined by the resolvent operator associated with $\mu \partial g$ and is single-valued~\cite{parboy13}, we have
	\be
	x^\star
    	\; - \;
   	 \prox_{\mu g}(x^\star \, -\, \mu \nabla f(x^\star))
	 \; = \; 
	 0.
	 \label{eq.nco}
	\ee
We next demonstrate that~\eqref{eq.main11} can be solved} using the proximal gradient flow dynamics, \tc{black}{$\dot{x} = - \mu \, G_{\mu}(x)$},
	\be
    \ba{rcl}
    \dot{x}
    & \, = \, &
    -
    \left(
    x
    \,-\,
    \prox_{\mu g}(x\, -\, \mu \nabla f(x))
    \right)
    \\[0.cm]
    & \, = \, &
    -
    \tc{black}{\mu}
    \left(\nabla f(x)\,+\, \nabla M_{\mu g} (x \,-\, \mu\nabla f(x)) \right).
    \ea
    \label{eq.pgGF}
    \ee

	\begin{remark}
Proximal gradient flow dynamics~\eqref{eq.pgGF} are different from the subgradient flow dynamics associated with nonsmooth problem~\eqref{eq.main11}. Standard proximal gradient algorithm~\cite{becteb09} can be obtained via explicit forward Euler discretization of~\eqref{eq.pgGF}  with \tc{black}{the stepsize one, $x^{k+1} = \prox_{\mu g}(x^k - \mu \nabla f(x^k))$.} This should be compared and contrasted with~\cite[Section~4.1.1]{parboy13} in which implicit backward Euler discretization of the subgradient flow dynamics associated with~\eqref{eq.main11} was used. We also note that~\eqref{eq.pgGF} can be obtained by substituting $-\nabla f(x)$ for the dual variable $y$ in the $x$-update step of primal-descent dual-ascent gradient flow dynamics~\eqref{eq.GF} with $T = I$.
	\end{remark}
	
\tc{black}{In Sections~\ref{sec.SC} and~\ref{sec.cvx},} we examine properties of system~\eqref{eq.pgGF}, first for strongly convex problems and then for the problems in which only the PL condition holds.

	\begin{figure}
	\centering
    	\hspace*{-3.25cm} {\scalebox{0.725}{%_______________________________________________________________________________
%
%   Block diagram of the periodic modification to the original dynamics
%   drawn from Right to Left
%
%   Mihailo Jovanovic, February 17, 2020
%_______________________________________________________________________________
%
% TikZ styles for drawing
%
\input{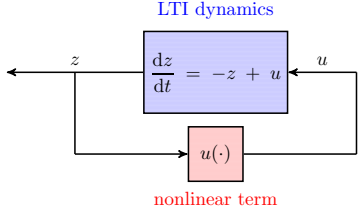}
%
%   set a filename for externalization
% \tikzsetnextfilename{clp_2dof_input_pert_config}
%
\noindent
\begin{tikzpicture}[scale=1, auto, >=stealth']

	% output node
	% starting point for uend
	% \node [input, name=uend] {};
	\node[] (end) at (0,0) {};
			
	% \node[] (midpoint1) at ($(uend) + (1.5cm,-1.75cm)$) {};
	 
   % original dynamics
    \node[block, minimum height = 1.5cm, top color=blue!20, bottom color=blue!20] (plant) at ($(end) + (6cm,0)$) {$ \dfrac{\mrd z}{\mrd t} \; = \; - z \; + \; u $};
    
    \node[] (uend) at ($(plant.west) - (1.25cm,0cm)$) {};
    
    \node[] (uend1) at ($(plant.west) - (2.5cm,0cm)$) {};
    
    % periodic feedback
    \node[block, minimum height = 1cm, top color=red!20, bottom color=red!20] (Gamma) at ($(end) + (6cm,-1.5cm)$) {$ u (\cdot)$};

	\node[] (ubegin) at ($(plant.east) + (1.25cm,0cm)$) {};
	
	\node[] (Gend) at ($(uend.center) - (0cm,1.5cm)$) {};

	\node[] (Gbegin) at ($(ubegin.center) - (0cm,1.5cm)$) {};

	% input nodes
	% inputs to plant
	% input dbegin to block plant
    	\draw [connector] (ubegin.center) -- node [midway, above] {$ u $} (plant.east);
		
	% connect plant with uenddown
         \draw [connector] (plant.west) -- node [midway, above] {$ z $} (uend1.center);
         
         % connect plant with uenddown
         \draw [line] (uend.center) -- (Gend.center);

        \draw [connector] (Gend.center) -- (Gamma.west);
                 
        \draw [line] (Gamma.east) -- (Gbegin.center);
        
        \draw [line] (Gbegin.center) -- (ubegin.center);
        
        	 \node [below = 0.05cm of Gamma](extra){\tc{red}{nonlinear term}}; 
	 
	  \node [above = 0.cm of plant](extra){\tc{blue}{$\ba{c} \mbox{LTI dynamics} \ea$}};  %

%         % connect plant with uenddown
%         \draw [line] (midpoint1.center) -- (midpoint2.center);
%
%	% connect midpoint1 with periodic block
%	 \draw [connector] (midpoint2.center) -- (periodic.west);
%	 
%	  \draw [line] (periodic.east) -- (midpoint3.center);
%
%         \draw [line] (midpoint3.center) -- (midpoint4.center);
%
%	 \draw [connector] (midpoint4.center) -- (plantdowneast.center);
%	 
%	 \node [below = -0.07cm of periodic](extra){\tc{red}{$\ba{c} \mbox{\bf sensor-less feedback} \ea$}};  %

\end{tikzpicture}
%_______________________________________________________________________________
}}
	\caption{\tc{black}{Both the proximal gradient flow dynamics~\eqref{eq.pgGF} and the DR splitting dynamics~\eqref{eq.DRSGF} can be represented via feedback interconnections of stable LTI systems with nonlinear terms that possess certain structural properties.}}
	\label{fig.bd}
	\end{figure}

	\vspace*{-2ex}
\subsection{Strongly convex problems}
\label{sec.SC}

	\vspace*{-2ex}
We utilize the theory of integral quadratic constraints to prove global asymptotic stability of the proximal gradient flow dynamics~\eqref{eq.pgGF}  under the following assumption.
	\vspace*{-2ex}
	\begin{assumption}
	\label{as.main}
Let $f$ in~\eqref{eq.main11} be $m_f$-strongly convex, let $\nabla f$ be $L_f$-Lipschitz continuous, and let the regularization function $g$ be proper, closed, and convex.
 	\end{assumption}
	
	\vspace*{-1ex}
\tc{black}{As illustrated in Fig.~\ref{fig.bd}, system~\eqref{eq.pgGF}} can be expressed as a feedback interconnection of an LTI system
	\begin{subequations}
	 \label{eq.sspgGF}
    \be
    \ba{rcl}
    \dot{z}
    & \, = \; &
    A\,z \; + \; B\,u,
    ~~
    \xi
    \; = \; 
    C \, z
   \\[0.cm]
    A
    & \, = \; &
    - I,
    ~~
    B \; = \; C \; = \; I
    \ea
    \label{eq.LTI}
    \ee
\tc{black}{where $z \DefinedAs x$,} with the nonlinear term,
	 \be
    u ( \xi )
    \; \DefinedAs \;
    \prox_{\mu g} ( \xi \,-\, \mu\,\nabla f (\xi) ).
    \label{eq.u-prox}
    \ee
Lemma~\ref{lem.K} combines firm nonexpansiveness of $\prox_{\mu g}$, strong convexity of $f$, and Lipschitz continuity of $\nabla f$ to characterize nonlinear map~\eqref{eq.u-prox} by establishing a quadratic inequality for
 	$
	u( \xi ).
	$
 	\end{subequations}
	\vspace*{-2ex}
\begin{lemma}
	\label{lem.K}
Let Assumption~\ref{as.main} hold. Then, for any $\xi \in \bbR^n$, $\hat{\xi} \in \bbR^n$, $
	u
    \DefinedAs
    \prox_{\mu g} ( \xi - \mu\,\nabla f(\xi) ),
	$
and $\hat{u} \DefinedAs \prox_{\mu g} (\hat{\xi} - \mu\,\nabla f( \hat{\xi} ))$, the pointwise quadratic inequality
	\begin{subequations}
	\label{eq.Pi}
	 \be
    \tbo{\xi \, - \, \hat{\xi}}{u \, - \, \hat{u}}^T
    \underbrace{\tbt{\sigma^2 I}{0}{0}{-I}}_{\Pi}
    \tbo{\xi  \, - \,  \hat{\xi}}{u  \, - \, \hat{u}}
    \, \geq \,
    0
    \label{eq.Pi}
    \ee
holds, where
    \be
    \sigma
    \;=\;
    \max \left\{|1\,-\,\mu m_f|,|1\,-\,\mu L_f|\right\}.
    \label{eq.tau}
    \ee
    \end{subequations}
\tc{black}{Moreover, the nonlinear function $u (\xi) \DefinedAs \prox_{\mu g} (\xi - \mu \nabla f(\xi))$ is a contraction for $\mu \in (0,2/L_f)$.}
\end{lemma}
	\vspace*{-4ex}

\begin{pf} \tc{black}{Since $\prox_{\mu g}$ is firmly nonexpansive~\cite{parboy13}, it is also Lipschitz continuous with parameter $1$, i.e., 
    \be
    \|
    u
    \, - \,
    \hat{u}
    \|^2_2
    \, \leq \,
    \|( \xi - \mu \nabla f(\xi))\,-\, (\hat{\xi} - \mu \nabla f(\hat{\xi}))\|_2^2.
    \label{eq.prox-ne}
    \ee
Expanding the right-hand-side of~\eqref{eq.prox-ne} yields,
    \be
    \ba{rcl}
    \| u\, - \,\hat{u}\|^2_2
    & \! \leq \! &
    \| \xi \, - \, \hat{\xi} \|_2^2
    \; + \;
    \mu^2 \| \nabla f ( \xi ) \, - \, \nabla f ( \hat{\xi} )\|_2^2
    ~ -
    \\[.1cm]
    & \!  \! &
    2 \mu \inner{ \xi \, - \, \hat{\xi}}{\nabla f(\xi)\,-\, \nabla f(\hat{\xi})}
    \ea
    \non
    \ee
and utilizing inequality~\eqref{eq.Lfmf} for an $m_f$-strongly convex function $f$ with an $L_f$-Lipschitz continuous gradient, the last inequality can be further simplified to obtain,
    \be
    \ba{rcl}
    \|u\, - \,\hat{u}\|^2_2
    & \leq &
    (1\,-\,\dfrac{2 \, \mu \, m_f L_f}{L_f + m_f}) \, \| \xi \, - \, \hat{\xi}\|_2^2
    ~+
    \\[.25cm]
    & &
    (\mu^2\,-\,\dfrac{2\,\mu}{L_f + m_f}) \, \|\nabla f(\xi) \, - \, \nabla f(\hat{\xi})\|_2^2.
     \ea
    \label{eq.tempmu}
    \ee
Depending on the sign of $\mu - {2}/{(L_f+m_f)}$ \tc{black}{either lower or upper bound in~\eqref{eq.Lf-mf}} can be used to upper bound the second term on the right-hand-side of~\eqref{eq.tempmu}, thereby yielding
     \be
    \|u\, - \,\hat{u}\|^2_2
    \, \leq \,
    \max
    \left\{
    (1\,-\,\mu L_f)^2,(1\,-\,\mu m_f)^2
    \right\}
    \| \xi \,-\,\hat{\xi} \|_2^2.
    \label{eq.tempmu1}
    \ee
Thus, for $\sigma$ given by~\eqref{eq.tau} the nonlinear function $u (\xi)$ is a contraction if and only if
    $
    -1  < 1 - \mu L_f < 1
    $
 and 
 	$  
     -1 <  1 - \mu m_f <  1.
    $
Since $m_f \leq L_f$, these conditions hold for $\mu \in ( 0, 2/L_f)$ which completes the proof.}
\end{pf}

	\vspace*{-4ex}
\tc{black}{We next employ~\cite[Theorem 3]{husei16} to prove the global exponential stability of the equilibrium point $z^\star$ of~\eqref{eq.sspgGF} with the rate $\rho > 0$ by verifying the existence of a positive definite matrix $P$ such that,
\be
	\tbt{A_\rho^TP + PA_\rho}{PB}{B^TP}{0}
	\, + \,
	\tbt{C^T}{0}{0}{I} 
	\Pi 
	\tbt{C}{0}{0}{I}
	\, \preceq \,
	0,
	\label{eq.KYP}
\ee
where $A_\rho \DefinedAs A + \rho I$ and $\Pi$ is given by~\eqref{eq.Pi}.}

\begin{theorem}
\label{thm.pg}
Let Assumption~\ref{as.main} hold and let $\mu \in (0,2/L_f)$. Then, \tc{black}{the equilibrium point $z^\star$ of the} proximal gradient flow dynamics~\eqref{eq.sspgGF} \tc{black}{is globally $\rho$-exponentially stable, i.e., there is $c > 0$ and $\rho \in ( 0, 1-\sigma ]$ such that,} 
	\[
	\| z (t) \, - \, z^\star \|_2
	 \; \leq \;
	c \, \mre^{-\rho t}
	\,
	\norm{z(0) \, - \, z^\star}_2,
	~~
	\forall \, t \, \geq \, 0
	\]
where $\sigma$ is given by~\eqref{eq.tau}. \tc{black}{Moreover, $x^\star = z^\star$ is the optimal solution of~\eqref{eq.main11}.}
\end{theorem}
	\vspace*{-4ex}
\begin{pf}
Substituting $\Pi$ given by~\eqref{eq.Pi} into~\eqref{eq.KYP} implies that the condition~\eqref{eq.KYP} holds if \tc{black}{there exists a positive scalar $p$ such that   
	\be
    \tbt{2 (1 - \rho) p - \sigma^2}{-p}{-p}{1}
    \; \succeq \;
    0
    \label{eq.temp4}
    \ee
where the block-diagonal structure of $A$, $B$, $C$, and $\Pi$ allows us to choose $P = p I$ without loss of generality.  Condition~\eqref{eq.temp4} is satisfied if there is $p > 0$ such that 
	\be
	p^2 
	\; - \; 
	2(1 \, - \, \rho) p
	\; + \;
	\sigma^2 
	\; \leq \;
	0
	\label{eq.p}
	\ee
where $\rho < 1$ guarantees positivity of the first element on the main diagonal of the matrix in~\eqref{eq.temp4}. For $\mu \in (0,2/L_f)$, Lemma~\ref{lem.K} implies $\sigma < 1$ and $\rho \leq 1 - \sigma$ is required for the existence of $p > 0$ such that~\eqref{eq.p} holds. Thus, $z^\star$ is globally exponentially stable with the rate $\rho \leq 1-\sigma$. The result follows because the equilibrium point $z^\star = x^\star$ of~\eqref{eq.sspgGF} satisfies the optimality condition~\eqref{eq.nco} for optimization problem~\eqref{eq.main11}.}
	\end{pf}

\begin{remark}
For $\mu =2/(L_f+m_f)$, the second term on the right-hand-side in~\eqref{eq.tempmu} disappears and $\sigma$ is given by
	$
    \sigma
    = 
    (L_f - m_f)/(L_f + m_f)
    = 
    (\kappa - 1)/(\kappa + 1)
    $
where $\kappa \DefinedAs L_f/m_f$ is the condition number of the function $f$ and $\rho$ is upper bounded by \tc{black}{$2/(\kappa + 1)$. In fact, this is the best achievable convergence rate for system~\eqref{eq.pgGF}.}
\end{remark}

	\vspace*{-2ex}
\subsection{Proximal Polyak-Lojasiewicz condition}
	\label{sec.cvx}

		\vspace*{-2ex}
Next, we consider the problems in which the function $f$ is not strongly convex but the function $ F \DefinedAs f + g $ satisfies the proximal PL condition~\eqref{eq.lbonG}.

\begin{assumption}
	\label{as.pl}
Let the regularization function $g$ in~\eqref{eq.main11} be proper, closed, and convex, let $f$ be twice continuously differentiable with $\nabla^2 f (x) \preceq L_f I$, and let the generalized gradient map satisfy the proximal PL condition,
	\be
    \|G_{\mu}(x)\|_2^2
    \,\geq\,
    \gamma\, (F_{\mu} (x)\,-\, F_{\mu}^\star)
    \non
    \ee
where $\mu \in (0,1/L_f)$, $\gamma>0$, and $F_{\mu}^{\star}$ is the optimal value of the FB envelope $F_{\mu}$.
\end{assumption}

\begin{remark}
The proximal gradient algorithm can be interpreted as a variable-metric gradient method on FB envelope and~\eqref{eq.pgGF} can be equivalently written as
    \[
    \dot{x}
    \;=\;
    - \tc{black}{\mu} \, (I\,-\,\mu \nabla^2 f(x))^{-1}\,\nabla F_{\mu}(x).
    \]
Under Assumption~\ref{as.pl}, 
	$
	I - \mu \nabla^2 f(x)
	$
is invertible and the functions $F$ and $F_{\mu}$ have the same minimizers and the same optimal values~\cite{patstebem14}, i.e.,
	$
	\argmin_x F (x)
          =
          \argmin_x F_{\mu} (x)
          $
and
   	$
	F^\star
         =
         F_{\mu}^\star.
    	$
This motivates the analysis of the convergence properties of~\eqref{eq.pgGF} in terms of the FB envelope.
\end{remark}
\begin{theorem}
Let Assumption~\ref{as.pl} hold. Then the forward-backward envelope associated with the proximal gradient flow dynamics~\eqref{eq.pgGF} converge exponentially to $F_\mu^\star = F^\star$ with the rate $\rho = \gamma \tc{black}{\mu} (1-\mu L_f)$,~i.e.,
	\[
	F_{\mu} (x(t)) \, - \, F_{\mu}^\star
	 \; \leq \;
	\mre^{-\rho t}
	\,
	(
	F_\mu ( x(0) ) \, - \, F_{\mu}^\star
	),
	~~
	\forall \, t \, \geq \, 0.
	\]	
\end{theorem}
		\vspace*{-6ex}
\begin{pf}
For a Lyapunov function candidate,
    \[
    V(x)
    \;=\;
    F_{\mu}(x)\,-\,F_{\mu}^\star
    \]
the derivative of $V$ along the solutions of~\eqref{eq.pgGF} is given by
    \be
    \ba{rcl}
    \dot{V}(x)
    & \, = \, &
    \inner{\nabla F_{\mu}(x)}{\dot{x}}
    \\[0.cm]
	& \, = \, &
     -\, \inner{\nabla F_{\mu}(x)}{\tc{black}{\mu} (I\,-\,\mu \nabla^2 f(x))^{-1} \nabla F_{\mu}(x)}
         \\[0.cm]
	& \, = \, &
     -\, \inner{G_{\mu}(x)}{\tc{black}{\mu} (I\,-\,\mu \nabla^2 f(x))\,G_{\mu}(x)}.
    \ea
    \label{eq.dLyap}
    \non
    \ee
Since the gradient of $f$ is $L_f$-Lipschitz continuous, i.e., $\nabla^2 f (x) \preceq L_f I$ for all $x \in \bbR^n$, Assumption~\ref{as.pl} implies
	$
	-(I - \mu \nabla^2 f(x))
	\preceq 
	- (1 - \mu L_f) I,
	$
and, thus, 	
	\be
	\ba{rcl}
    \dot{V}(x)
   & \, \leq \, &
    -
    \tc{black}{\mu} (1\,-\,\mu L_f) \,\|G_{\mu}(x)\|_2^2
    \\[0.cm]
    & \, \leq \, &
    -
    \gamma
    \tc{black}{\mu}
    (1\,-\,\mu L_f)
    \,
    (F_{\mu} (x)\,-\, F_{\mu}^\star)
    \ea
    \label{eq.temp2}
    \ee
is non-positive for $\mu \in (0,1/L_f)$. Moreover, combining the last inequality with the definition of $V$ yields
	$
    \dot{V}
      \leq 
    -
    \gamma
    \tc{black}{\mu} (1 - \mu L_f)
    V,
    $
 which implies
	\be
	F_{\mu} (x (t))\,-\,F_{\mu}^\star
	\; \leq \;
	\mre^{-
    	\gamma
    	\tc{black}{\mu}
    	(1-\mu L_f) t}
	\,
	(F_{\mu}(x (0))\,-\,F_{\mu}^\star).
	\non
	\ee
\end{pf}

\begin{remark}
When the proximal PL condition is satisfied, $F_{\mu}(x (t))- F_{\mu}^{\star}$ converges exponentially but, in the absence of strong convexity, the exponential convergence rate cannot be established for $\norm{x(t) - x^{\star}}_2$. Thus, although the objective function converges exponentially fast, the solution to~\eqref{eq.pgGF} does not enjoy this convergence rate. To the best of our knowledge, the convergence rate of $x(t)$ to the set of optimal values $x^\star$ is not known in this case.
\end{remark}

 \vspace*{-3ex}
\section{Global exponential stability of the Douglas-Rachford splitting dynamics}
\label{section.DRS}

	 \vspace*{-2ex}
We next introduce a continuous-time gradient flow dynamics based on the well-known Douglas-Rachford splitting  algorithm~\cite{dourac56} and establish global exponential stability for strongly convex $f$.

\vspace*{-3ex}
\subsection{Non-smooth composite optimization problem}

	\vspace*{-2ex}
The optimality condition for~\eqref{eq.main11} is given by~\eqref{eq.OC}, i.e., 
   $
    0
    \in
    \nabla f (x^\star) + \partial g(x^\star).
    $
Multiplication by $\mu$ and \tc{black}{addition of $x$ to the both sides} yields
	$
       	0
       	\in 
    	\left[
    	I + \mu \nabla f
   	 \right] (x^\star)
    	+ 
    	\mu \partial g(x^\star) - x^\star.
    $
Since $\prox_{\mu f} \DefinedAs (I + \mu \nabla f)^{-1}$ \tc{black}{is single-valued~\cite{parboy13},} introducing $z \DefinedAs x- \mu \partial g(x)$ leads to,
\begin{subequations}
\label{xz-opt}
    \be
    x^\star
    \; = \;
    \prox_{\mu f}(x^\star \, - \, \mu \partial g(x^\star))
    \; = \;
    \prox_{\mu f}(z^\star).
    \label{eq.xz-opt1}
    \ee
Now, adding $x$ to the both sides of the defining equation for $z$ gives
    $
    \left[ I + \mu \partial g \right] (x^\star)
    = 
    2 \prox_{\mu f}(z^\star) - z^\star,
   $
i.e., 
    \be
    x^{\star} \;=\; \prox_{\mu g} (2\,\prox_{\mu f}(z^{\star})\,-\,z^{\star}).
    \label{eq.xz-opt2}
    \ee
Combining~\eqref{eq.xz-opt1} and~\eqref{eq.xz-opt2} results in the following optimality condition,     
  	\be
    	\prox_{\mu f}(z^{\star})
	 \; - \;
	 \prox_{\mu g} (2\,\prox_{\mu f}(z^{\star})\,-\,z^{\star})
	 \; = \; 
	 0.
    \label{eq.xz}
    \ee
Furthermore, the reflected proximal operators~\cite{gisboy17},
	 $
    R_{\mu f} (z)
    \DefinedAs 
    [ 2\,\prox_{\mu f} - I ] (z)
    $
    and
    $
    R_{\mu g}
     \DefinedAs 
    [ 2\,\prox_{\mu g} - I ] (z),
    $
can be used to rewrite optimality \mbox{condition~\eqref{eq.xz} as}
	\be
	z^\star
	\; - \; 
	[ R_{\mu g} R_{\mu f} ] (z^\star)
	\; = \;
	0.
	\label{eq.xz1}
	\ee
	\end{subequations}

\tc{black}{We are now ready to introduce the continuous-time DR gradient flow dynamics to compute $z^\star$,
    \be
    % \ba{rcl}
    \dot{z}
    \; = \;
    -z \; + \; [ R_{\mu g} R_{\mu f} ] (z).
    \label{eq.DRSGF}
    \ee
Note that the explicit forward Euler discretization of~\eqref{eq.DRSGF} yields the standard DR splitting algorithm~\cite{eckber92}.}
\tc{black}{We view~\eqref{eq.DRSGF}} as a feedback interconnection of an \tc{black}{LTI system~\eqref{eq.LTI}} with the nonlinear term,
	\be
    u(\xi)
    \; \DefinedAs \;
    [ R_{\mu g} R_{\mu f} ] (\xi).
    \label{eq.u-z}
    \ee
We first characterize properties of nonlinearity $u$ in~\eqref{eq.u-z} and then, similar to the previous section, establish global exponential stability of \tc{black}{nonlinear system~\eqref{eq.DRSGF}.}

\begin{lemma}
Let Assumption~\ref{as.main} hold and let $\mu \in (0,2/L_f)$. Then, the operator $R_{\mu f}$ is $\sigma$-contractive,
    \[
    \|R_{\mu f}(x) \,-\,R_{\mu f}(y)\|_2
    \,\leq\,
    \sigma
    \|x\,-\,y\|_2
    \]
where $\sigma$ is given by
     \be
     \sigma \,=\, \max \left\{|1\,-\,\mu m_f|,|1\,-\,\mu L_f|\right\}\,<\,1.
     \label{eq.sigma}
     \ee
\label{thm.contract}
\end{lemma}
		\vspace*{-6ex}
\begin{pf}
Given $z_x \DefinedAs \prox_{\mu f}(x)$ and $z_y \DefinedAs \prox_{\mu f}(y)$, $x$ and $y$ can be computed as follows
    \[
    x\;=\; z_x\,+\,\mu\,\nabla f(z_x),
    ~
    y\;=\; z_y\,+\,\mu\,\nabla f(z_y).
    \]
Thus,
    \[
    \ba{l}
    \|R_{\mu f}(x)\,-\,R_{\mu f}(y)\|^2
    \,=\,
    \|2(z_x - z_y) \,-\, (x  - y)\|^2\; =
    \\[0.cm]
    \|(z_x - z_y) \,-\, \mu\,(\nabla f(z_x) - \nabla f(z_y))\|^2
    \; = \,
    \|z_x - z_y\|^2\, +\,
    \\[0.cm]
    \tc{black}{\|}\mu (\nabla f(z_x) - \nabla f(z_y))\|^2
    -
    2\,\mu\inner{\nabla f(z_x) - \nabla f(z_y)}{z_x - z_y}
    \\[0.cm]
    \leq \,
    \max
    \left\{
    (1\,-\,\mu L_f)^2,(1\,-\,\mu m_f)^2
    \right\}
    \,\|z_x\,-\,z_y\|^2
    \\[0.cm]
    \leq \,
    \tc{black}{\sigma^2}
    \,\|x\,-\,y\|^2.
    \ea
    \]
where the firm non-expansiveness of $\prox_{\mu f}$ is used in the last step. Moreover, according to Lemma~\ref{lem.K}, for $\mu \in (0,2/L_f)$ we have $\sigma<1$, which completes the proof.
\end{pf}

\begin{lemma}
Let Assumption~\ref{as.main} hold and let $\mu \in (0,2/L_f)$. Then, the operator $R_{\mu g}$ is firmly non-expansive.
\end{lemma}
		\vspace*{-4ex}
\begin{pf}
	\\
	$
        \|R_{\mu g}(x) - R_{\mu g}(y)\|_2^2 
        =
    4 \|\prox_{\mu f}(x) - \prox_{\mu f}(y)\|_2^2
    +
    \|x - y\|_2^2
    -
    4\inner{x - y}{\prox_{\mu f}(x) - \prox_{\mu f}(y)}
    \leq
    \|x - y\|_2^2.
    $
   \end{pf}

\begin{remark}
\label{re.IQCs}
Since $R_{\mu g}$ is firmly non-expansive and $R_{\mu f}$ is $\sigma$-contractive, the composite operator $R_{\mu g} R_{\mu f}$ is also $\sigma$-contractive. Moreover, since the operator $R_{\mu f}$ and nonlinearity $u$ in~\eqref{eq.u-prox} have the same contraction parameters, the quadratic inequality that describes~\eqref{eq.u-prox} can be also used to characterize the composite operator $R_{\mu g} R_{\mu f}$.
\end{remark}

\begin{theorem}
\label{thm.DRS}
Let Assumption~\ref{as.main} hold \tc{black}{and let $\mu \in (0,2/L_f)$.} Then, \tc{black}{the equilibrium point $z^\star$ of} the DR splitting dynamics~\eqref{eq.DRSGF} \tc{black}{is} globally \tc{black}{$\rho$-}exponentially stable, i.e., there is $c > 0$ and $\rho \in (0, 1-\sigma)$ such that,
	\[
	\norm{z (t) \, - \, z^\star}
	 \; \leq \;
	c \, \mre^{-\rho t} \,
	\norm{z(0) \, - \, z^\star},
	~~
	\forall \, t \, \geq \, 0
	\]
\tc{black}{where $\sigma$ is given by~\eqref{eq.sigma}. Moreover, $x^\star = \prox_{\mu f}(z^\star)$ is the optimal solution of~\eqref{eq.main11}.}
\end{theorem}
	\vspace*{-4ex}
\begin{pf}
Although the nonlinear terms in systems~\eqref{eq.GF} and~\eqref{eq.DRSGF} are different, \tc{black}{they share quadratic characterization~\eqref{eq.Pi} and the LTI dynamics~\eqref{eq.LTI}. Thus, the result follows from the proof of Theorem~\ref{thm.pg} and the fact that $x^\star = \prox_{\mu f}(z^\star)$ satisfies optimality condition~\eqref{eq.xz}.}
	\end{pf}

\vspace*{-2ex}
\subsection{Douglas-Rachford splitting on the dual problem}

	\vspace*{-2ex}
\tc{black}{Even though} the DR splitting algorithm cannot be directly used to solve a problem with a more general linear equality constraint,
    \be
    \ba{rcl}
	&&\minimize\limits_{x, \, z}
	~~~~
	f(x)
	\; + \;
	g(z)
    \\[.1cm]
    &&\subject
    ~~~
    T x\,+\,S z\,=\,r
    \ea
    \label{eq.genproblast}
	\ee
it can be utilized to solve the dual problem,
    \beq
    \ba{rcl}
	&&\minimize\limits_{\zeta}
	~~~
	f_1(\zeta)
	\, + \,
	g_1(\zeta).
    \ea
    \label{eq.dualgen1}
	\eeq
\tc{black}{Here, $T \in \bbR^{m \times n}$, $S \in \bbR^{m \times n}$, and $r \in \bbR^{m}$ are the problem parameters,
	$
	f_1(\zeta)
        \DefinedAs
    	f^{\star}(-T^{T}\zeta)
        +
    	r^T \zeta,
	$
	$
    	g_1(\zeta)
    	\DefinedAs
	g^{\star}(-S^{T}\zeta),
   	$}
and $h^{\star}(\zeta)\DefinedAs \sup_x ( \zeta^T x-h(x) )$ is the conjugate of the function $h$. It is a standard fact~\cite{eckber92,gab83} that solving the dual problem~\eqref{eq.dualgen1} via the DR splitting algorithm is equivalent to using ADMM for the original problem~\eqref{eq.genproblast}. \tc{black}{If Assumption~\ref{as.main} holds and if $T$ is a full row rank matrix, the global exponentially stability of the DR gradient flow dynamics associated with~\eqref{eq.dualgen1}, 
	$
	\dot{\zeta}
    	= 
    	- \zeta + [ R_{\mu g_1} R_{\mu f_1} ] (\zeta),
	$
is readily established.}

	\vspace*{-2ex}
\section{Concluding remarks}
	\label{sec.conclusion}

	\vspace*{-2ex}
We study a class of nonsmooth optimization problems in which it is desired to minimize the sum of a continuously differentiable function with a Lipschitz continuous gradient and a nondifferentiable function. For strongly convex problems, we employ the theory of integral quadratic constraints to prove global exponential stability of proximal gradient flow and Douglas-Rachford splitting dynamics. We also utilize a generalized Polyak-Lojasiewicz condition for nonsmooth problems to demonstrate the global exponential convergence of the forward-backward envelope for the proximal gradient flow algorithm even in the absence of strong convexity.

\appendix

	\vspace*{-2ex}
\section{Proximal PL condition}
	\label{sec.lbonG}

		\vspace*{-2ex}
The generalization of the PL condition to nonsmooth problems was introduced in~\cite{karnutsch16} and is given by
    \be
    \mathcal{D}_g(x,L_f)
    \;\geq\;
    2 \kappa\,(F(x) \,-\, F^{\star})
    \label{eq.PPL}
    \ee
where $\kappa$ is a positive constant, $L_f$ is the Lipschitz constant of $\nabla f$, and $\mathcal{D}_g(x,\alpha)$ is determined by
    \be
    -2\alpha\min\limits_{y}
    \,
    (
    \inner{\nabla f(x)}{y - x}
    +
    \dfrac{\alpha}{2}\,\|y - x\|_2^2
    +
    g(y) - g(x)
    ).
    \label{eq.Dg}
    \ee
Herein, we show that if proximal PL condition~\eqref{eq.PPL} holds, there is a lower bound given by~\eqref{eq.lbonG} on the norm of the generalized gradient map $G_{\mu}(x)$. For $\mu \in (0,1/L_f)$,
	$
	\tc{black}{\mathcal{D}_g}(x,1/\mu)
         \geq
    	\tc{black}{\mathcal{D}_g}(x,L_f),
   	$
and, thus, inequality~\eqref{eq.PPL} also holds for $\tc{black}{\mathcal{D}_g}(x,1/\mu)$. Moreover, from the definition~\eqref{eq.Dg} of $\tc{black}{\mathcal{D}_g}(x,\alpha)$, it follows that
    \[
    \tc{black}{\mathcal{D}_g}(x,1/\mu)
    \;=\;
    \dfrac{2}{\mu}\,(F(x)\,-\,F_{\mu}(x))
    \]
where $F \DefinedAs f + g$ and $F_{\mu}$ is the FB envelope. Substituting this expression for $\tc{black}{\mathcal{D}_g}(x,1/\mu)$ to~\eqref{eq.PPL} yields,
    \be
    \dfrac{1}{\mu}\,(F(x)\,-\,F_{\mu}(x))
    \;\geq\;
    \kappa\,(F(x)\,-\,F^\star).
    \label{eq.interm1}
    \ee
The smooth part of the objective function $f$ can be written~as~\cite{patstebem14},
    \be
    \ba{rcl}
    f(x)
    & \, = \;\, &
    F_{\mu}(x)
    \; - \;
    g(\prox_{\mu g}(x\,-\,\mu\,\nabla f(x)))
    ~+
    \\[0.1cm]
    & \, \;\, &
    \mu \inner{\nabla f(x)}{G_{\mu}(x)}
    \; - \;
    \dfrac{\mu}{2}\,\|G_{\mu}(x)\|_2^2
    \ea
    \non
    \ee
and substituting this expression for $f$ to~\eqref{eq.interm1} yields
    \be
    \ba{c}
    \tfrac{(\mu\kappa\,-\,1)}{2}\,\|G_{\mu}(x)\|_2^2
    \, \geq \,
    \kappa\,(F_{\mu}(x) - F^{\star})
    \,+\,
    \tfrac{(\mu\kappa - 1)}{\mu}\,g(x)
    ~ -
    \\[0.1cm]
    (\mu\kappa - 1) 
   (
    \tfrac{1}{\mu}\,g(\prox_{\mu g}(x - \mu\,\nabla f(x)))
    + 
    \inner{\nabla f(x)}{G_{\mu}(x)}
    ).
    \ea
    \label{eq.temp14}
    \ee
Since
    $
    G_{\mu}(x)
    -
    \nabla f(x)
    \in
    \partial g(x),
    $
the subgradient inequality~\eqref{eq.subgrad} implies
    \be
    \ba{rcl}
    0
    \;\leq\;
    \mu \,\|G_{\mu}(x)\|_2^2
    & \! \leq \! &
    g(x)
    \,-\,
    g(\prox_{\mu g}(x\,-\,\mu\,\nabla f(x)))
    ~ +
    \\[0.cm]
    & \!\! &
    \mu \inner{\nabla f(x)}{G_{\mu}(x)}.
    \ea
    \label{eq.temp15}
    \ee
\tc{black}{Combining~\eqref{eq.temp14} and~\eqref{eq.temp15} and taking the sign of $\mu\kappa - 1$ into account yields,
    \[
    \dfrac{\alpha}{2} \, \|G_{\mu}(x)\|_2^2
    \; \geq \;
    \kappa\,(F_{\mu}(x)\,-\,F^{\star}),
    ~~
    \alpha 
    \; \DefinedAs \; 
    | \mu \kappa \, - \, 1 |.
    \]}
Furthermore, since~\cite{patstebem14},
	$
	\argmin_x F (x)
          =
          \argmin_x F_{\mu} (x)
          $
and
   	$
	F^\star
         =
         F_{\mu}^\star,
    	$	
	$F_{\mu}^{\star}$ can be substituted for $F^{\star}$ and we have
    $
    \|G_{\mu}(x)\|_2^2
    \geq
    \gamma\, (F_{\mu} (x)\,-\, F_{\mu}^\star)
    $
with $\gamma \DefinedAs 2\kappa/| \mu \kappa - 1 |$.

\vspace*{-2ex}


\begin{thebibliography}{10}
\providecommand{\url}[1]{#1}
\csname url@rmstyle\endcsname
\providecommand{\newblock}{\relax}
\providecommand{\bibinfo}[2]{#2}
\providecommand\BIBentrySTDinterwordspacing{\spaceskip=0pt\relax}
\providecommand\BIBentryALTinterwordstretchfactor{4}
\providecommand\BIBentryALTinterwordspacing{\spaceskip=\fontdimen2\font plus
\BIBentryALTinterwordstretchfactor\fontdimen3\font minus
  \fontdimen4\font\relax}
\providecommand\BIBforeignlanguage[2]{{%
\expandafter\ifx\csname l@#1\endcsname\relax
\typeout{** WARNING: IEEEtran.bst: No hyphenation pattern has been}%
\typeout{** loaded for the language `#1'. Using the pattern for}%
\typeout{** the default language instead.}%
\else
\language=\csname l@#1\endcsname
\fi
#2}}

\bibitem{nedozd09}
A.~Nedi\'c and A.~Ozdaglar, ``Distributed subgradient methods for multiagent
  optimization,'' \emph{IEEE Trans. on Automat. Control}, vol.~54, no.~1, pp.
  48--61, 2009.

\bibitem{waneli11}
J.~Wang and N.~Elia, ``A control perspective for centralized and distributed
  convex optimization,'' in \emph{Proceedings of the 50th IEEE Conference on
  Decision and Control}, 2011, pp. 3800--3805.

\bibitem{latstepat16}
P.~Latafat, L.~Stella, and P.~Patrinos, ``New primal-dual proximal algorithm
  for distributed optimization,'' in \emph{Proceedings of the 55th IEEE
  Conference on Decision and Control}, 2016, pp. 1959--1964.

\bibitem{latfrepat19}
P.~Latafat, N.~Freris, and P.~Patrinos, ``A new randomized block-coordinate
  primal-dual proximal algorithm for distributed optimization,'' \emph{IEEE
  Trans. Automat. Control}, 2019, doi:10.1109/TAC.2019.2906924.

\bibitem{becteb09}
A.~Beck and M.~Teboulle, ``A fast iterative shrinkage-thresholding algorithm
  for linear inverse problems,'' \emph{SIAM J. Imaging Sci.}, vol.~2, no.~1,
  pp. 183--202, 2009.

\bibitem{parboy13}
N.~Parikh and S.~Boyd, ``Proximal algorithms,'' \emph{Found. Trends Opt.},
  vol.~1, no.~3, pp. 123--231, 2013.

\bibitem{arrhuruza58}
K.~J. Arrow, L.~Hurwicz, and H.~Uzawa, ``Studies in linear and non-linear
  programming,'' 1958.

\bibitem{feipagAUT10}
D.~Feijer and F.~Paganini, ``Stability of primal--dual gradient dynamics and
  applications to network optimization,'' \emph{Automatica}, vol.~46, no.~12,
  pp. 1974--1981, 2010.

\bibitem{chemalcorSCL16}
A.~Cherukuri, E.~Mallada, and J.~Cort{\'e}s, ``Asymptotic convergence of
  constrained primal--dual dynamics,'' \emph{Syst. Control Lett.}, vol.~87, pp.
  10--15, 2016.

\bibitem{chemallowcor18}
A.~Cherukuri, E.~Mallada, S.~Low, and J.~Cortes, ``The role of convexity on
  saddle-point dynamics: {L}yapunov function and robustness,'' \emph{IEEE
  Trans. Automat. Control}, vol.~63, no.~8, pp. 2449--2464, 2018.

\bibitem{dhikhojovTAC19}
N.~K. Dhingra, S.~Z. Khong, and M.~R. Jovanovi\'c, ``The proximal augmented
  {L}agrangian method for nonsmooth composite optimization,'' \emph{IEEE Trans.
  Automat. Control}, vol.~64, no.~7, pp. 2861--2868, July 2019.

\bibitem{quli18}
G.~Qu and N.~Li, ``On the exponential stability of primal-dual gradient
  dynamics,'' \emph{IEEE Control Syst. Lett.}, vol.~3, no.~1, pp. 43--48, 2018.

\bibitem{suboycan16}
W.~Su, S.~Boyd, and E.~Candes, ``A differential equation for modeling
  {N}esterov's accelerated gradient method: Theory and insights,'' \emph{J.
  Mach. Learn. Res.}, vol.~17, pp. 1--43, 2016.

\bibitem{wibwiljor16}
A.~Wibisono, A.~C. Wilson, and M.~I. Jordan, ``A variational perspective on
  accelerated methods in optimization,'' \emph{Proc. Natl. Acad. Sci.}, vol.
  113, no.~47, pp. E7351--E7358, 2016.

\bibitem{frarobvid18}
G.~Fran{\c{c}}a, D.~Robinson, and R.~Vidal, ``{ADMM} and accelerated {ADMM} as
  continuous dynamical systems,'' 2018, arXiv:1805.06579.

\bibitem{shidujorsu18}
B.~Shi, S.~Du, M.~I. Jordan, and W.~Su, ``Understanding the acceleration
  phenomenon via high-resolution differential equations,'' 2018,
  arXiv:1810.08907.

\bibitem{muejor19}
M.~Muehlebach and M.~I. Jordan, ``A dynamical systems perspective on {N}esterov
  acceleration,'' 2019, arXiv:1905.07436.

\bibitem{povli19}
J.~I. Poveda and N.~Li, ``Inducing uniform asymptotic stability in time-varying
  accelerated optimization dynamics via hybrid regularization,'' 2019,
  arXiv:1905.12110.

\bibitem{brobar89}
A.~Brown and M.~Bartholomew-Biggs, ``Some effective methods for unconstrained
  optimization based on the solution of systems of ordinary differential
  equations,'' \emph{J. Optimiz. Theory App.}, vol.~62, no.~2, pp. 211--224,
  1989.

\bibitem{schsin00}
J.~Schropp and I.~Singer, ``A dynamical systems approach to constrained
  minimization,'' \emph{Numer. Func. Anal. Opt.}, vol.~21, no. 3-4, pp.
  537--551, 2000.

\bibitem{zhamoksrajad18}
J.~Zhang, A.~Mokhtari, S.~Sra, and A.~Jadbabaie, ``Direct {R}unge-{K}utta
  discretization achieves acceleration,'' in \emph{Advances in Neural
  Information Processing Systems}, 2018, pp. 3900--3909.

\bibitem{lesrecpac16}
L.~Lessard, B.~Recht, and A.~Packard, ``Analysis and design of optimization
  algorithms via integral quadratic constraints,'' \emph{SIAM J. Optim.},
  vol.~26, no.~1, pp. 57--95, 2016.

\bibitem{huseiran17}
B.~Hu, P.~Seiler, and A.~Rantzer, ``A unified analysis of stochastic
  optimization methods using jump system theory and quadratic constraints,'' in
  \emph{Proceedings of the 2017 Conference on Learning Theory}, 2017, pp.
  1157--1189.

\bibitem{hules17}
B.~Hu and L.~Lessard, ``Dissipativity theory for {N}esterov's accelerated
  method,'' in \emph{Proceedings of the 34th International Conference on
  Machine Learning}, 2017, pp. 1549--1557.

\bibitem{fazribmorpre18}
M.~Fazlyab, A.~Ribeiro, M.~Morari, and V.~M. Preciado, ``Analysis of
  optimization algorithms via integral quadratic constraints: Nonstrongly
  convex problems,'' \emph{{SIAM} {J}. Optim.}, vol.~28, no.~3, pp. 2654--2689,
  2018.

\bibitem{mogjovACC18}
S.~Hassan-Moghaddam and M.~R. Jovanovi\'c, ``Distributed proximal augmented
  {L}agrangian method for nonsmooth composite optimization,'' in
  \emph{Proceedings of the 2018 American Control Conference}, Milwaukee, WI,
  2018, pp. 2047--2052.

\bibitem{mogjovCDC18a}
S.~Hassan-Moghaddam and M.~R. Jovanovi\'c, ``On the exponential convergence
  rate of proximal gradient flow algorithms,'' in \emph{Proceedings of the 57th
  IEEE Conference on Decision and Control}, Miami, FL, 2018, pp. 4246--4251.

\bibitem{dinhudhijovCDC18}
D.~Ding, B.~Hu, N.~K. Dhingra, and M.~R. Jovanovi\'c, ``An exponentially
  convergent primal-dual algorithm for nonsmooth composite minimization,'' in
  \emph{Proceedings of the 57th IEEE Conference on Decision and Control},
  Miami, FL, 2018, pp. 4927--4932.

\bibitem{seifazprepap19}
J.~Seidman, M.~Fazlyab, V.~Preciado, and G.~Pappas, ``A \mbox{control-theoretic
  approach to analysis and parameter selection} of {D}ouglas-{R}achford
  splitting,'' 2019, arXiv:1903.11525.

\bibitem{mohrazjovCDC18}
H.~Mohammadi, M.~Razaviyayn, and M.~R. Jovanovi\'c, ``Variance amplification of
  accelerated first-order algorithms for strongly convex quadratic optimization
  problems,'' in \emph{Proceedings of the 57th IEEE Conference on Decision and
  Control}, Miami, FL, 2018, pp. 5753--5758.

\bibitem{mohrazjovACC19}
H.~Mohammadi, M.~Razaviyayn, and M.~R. Jovanovi\'c, ``Performance of noisy
  {N}esterov's accelerated method for strongly convex optimization problems,''
  in \emph{Proceedings of the 2019 American Control Conference}, Philadelphia,
  PA, 2019, pp. 3426--3431.

\bibitem{mohrazjovTAC19}
H.~Mohammadi, M.~Razaviyayn, and M.~R. Jovanovi\'c, ``Robustness of accelerated
  first-order algorithms for strongly convex optimization problems,''
  \emph{IEEE Trans. Automat. Control}, 2019, submitted; also arXiv:1905.11011.

\bibitem{micschebe19}
S.~Michalowsky, C.~Scherer, and C.~Ebenbauer, ``Robust and structure exploiting
  optimization algorithms: An integral quadratic constraint approach,'' 2019,
  arXiv:1905.00279.

\bibitem{megran97}
A.~Megretski and A.~Rantzer, ``System analysis via integral quadratic
  constraints,'' \emph{IEEE Trans. Autom. Control}, vol.~42, no.~6, pp.
  819--830, 1997.

\bibitem{pol63}
B.~T. Polyak, ``Gradient methods for minimizing functionals,'' \emph{Zhurnal
  Vychislitel'noi Matematiki i Matematicheskoi Fiziki}, vol.~3, no.~4, pp.
  643--653, 1963.

\bibitem{patstebem14}
P.~Patrinos, L.~Stella, and A.~Bemporad, ``Forward-backward truncated {N}ewton
  methods for convex composite optimization,'' 2014, arXiv:1402.6655.

\bibitem{stethepat17}
L.~Stella, A.~Themelis, and P.~Patrinos, ``Forward--backward quasi-{N}ewton
  methods for nonsmooth optimization problems,'' \emph{Comput. Optim. Appl.},
  vol.~67, no.~3, pp. 443--487, 2017.

\bibitem{thestepat18}
A.~Themelis, L.~Stella, and P.~Patrinos, ``Forward-backward envelope for the
  sum of two nonconvex functions: Further properties and nonmonotone
  line-search algorithms,'' \emph{SIAM J. Optim.}, vol.~28, no.~3, pp.
  2274--2303, 2018.

\bibitem{dhikhojovTAC17}
N.~K. Dhingra, S.~Z. Khong, and M.~R. Jovanovi\'c, ``A second order primal-dual
  method for nonsmooth convex composite optimization,'' \emph{IEEE Trans.
  Automat. Control}, 2017, submitted; also arXiv:1709.01610.

\bibitem{nes13}
Y.~Nesterov, \emph{Introductory lectures on convex optimization: A basic
  course}, 2013, vol.~87.

\bibitem{karnutsch16}
H.~Karimi, J.~Nutini, and M.~Schmidt, ``Linear convergence of gradient and
  proximal-gradient methods under the {P}olyak-{L}ojasiewicz condition,'' in
  \emph{Joint European Conference on Machine Learning and Knowledge Discovery
  in Databases}, 2016, pp. 795--811.

\bibitem{husei16}
B.~Hu and P.~Seiler, ``Exponential decay rate conditions for uncertain linear
  systems using integral quadratic constraints,'' \emph{IEEE Trans. Autom.
  Control}, vol.~61, no.~11, pp. 3631--3637, 2016.

\bibitem{dourac56}
J.~Douglas and H.~Rachford, ``On the numerical solution of heat conduction
  problems in two and three space variables,'' \emph{{T}rans. {A}mer. {M}ath.
  {S}oc.}, vol.~82, no.~2, pp. 421--439, 1956.

\bibitem{gisboy17}
P.~Giselsson and S.~Boyd, ``Linear convergence and metric selection for
  {D}ouglas-{R}achford splitting and {ADMM},'' \emph{IEEE Trans. Automat.
  Control}, vol.~62, no.~2, pp. 532--544, 2017.

\bibitem{eckber92}
J.~Eckstein and D.~P. Bertsekas, ``On the {D}ouglas-{R}achford splitting method
  and the proximal point algorithm for maximal monotone operators,''
  \emph{Math. Program.}, vol.~55, no. 1-3, pp. 293--318, 1992.

\bibitem{gab83}
D.~Gabay, ``Applications of the method of multipliers to variational
  inequalities,'' in \emph{Studies in Mathematics and its Applications}, 1983,
  vol.~15, pp. 299--331.

\end{thebibliography}
\end{document}